\input amstex
\input amsppt.sty   
\hsize 30pc
\vsize 47pc
\def\nmb#1#2{#2}         
\def\cit#1#2{\ifx#1!\cite{#2}\else#2\fi} 
\def\idx{}               
\def\ign#1{}             

\redefine\o{\circ}

\predefine\ii\i
\redefine\i{^{-1}}
\define\row#1#2#3{#1_{#2},\ldots,#1_{#3}}

\def\today{\ifcase\month\or
 January\or February\or March\or April\or May\or June\or
 July\or August\or September\or October\or November\or December\fi
 \space\number\day, \number\year}
\topmatter
\title Characterizing algebras of smooth functions on manifolds
\endtitle
\author  Peter W. Michor \\
Ji\v r\'\ii{} Van\v zura
 \endauthor
\leftheadtext{\smc P\. Michor, J\. Van\v zura}
\affil
Erwin Schr\"odinger International Institute of Mathematical Physics,  
Wien, Austria \\
J\. Van\v zura: Mathematical Institute of the AV\v CR, department Brno,
Zi\v zkova 22, CZ 616 62 Brno, Czech Republic
\endaffil
\address 
P\. Michor: Institut f\"ur Mathematik, Universit\"at Wien, 
Strudlhofgasse 4, A-1090 Wien, Austria; and:  
Erwin Schr\"odinger International Institute of Mathematical Physics,  
Pasteurgasse 6/7, A-1090 Wien, Austria 
\endaddress 
\email Peter.Michor\@esi.ac.at \endemail
\address
J\. Van\v zura: Mathematical Institute of the AV\v CR, department Brno,
Zi\v zkova 22, CZ 616 62 Brno, Czech Republic
\endaddress
\email vanzura\@ipm.cz \endemail
\date {April 8, 1994} \enddate
\keywords $C^\infty$-algebra, smooth manifold\endkeywords
\subjclass 46J20, 51K10, 58A03, 58A05\endsubjclass
\abstract Among all $C^\infty$-algebras we characterize those which 
are algebras of smooth functions on smooth separable Hausdorff 
manifolds.
\endabstract
\endtopmatter

\document


\subhead\nmb0{1}. $C^\infty$-algebras \endsubhead
An $\Bbb R$-algebra is a commutative ring $A$ with unit together with 
a ring homomorphism $\Bbb R\to A$. Then every map 
$p:\Bbb R^n\to \Bbb R^m$ which is given by an $m$-tuple of real 
polynomials $(p_1,\dots,p_m)$ can be interpreted as a mapping 
$A(p):A^n\to A^m$ in such a way that projections, composition, and 
identity are preserved, by just evaluating each polynomial $p_i$ on 
an $n$-tuple $(a_1,\dots,a_n)\in A^n$.

A $C^\infty$-algebra $A$ is a real algebra in which we can moreover 
interpret all smooth mappings $f:\Bbb R^n\to \Bbb R^m$. There is a 
corresponding map $A(f):A^n\to A^m$, and again projections, 
composition, and the identity mapping are preserved. 

More precisely, a $C^\infty$-algebra $A$ is a product preserving 
functor from the category $C^\infty$ to the category of sets, where 
$C^\infty$ has as objects all spaces $\Bbb R^n$, $n\ge 0$, and all 
smooth mappings between them as arrows. Morphisms between 
$C^\infty$-algebras are then natural transformations: they correspond 
to those algebra homomorphisms which preserve the intepretation of 
smooth mappings.

This definition of $C^\infty$-algebras is due to Lawvere \cit!{2}, 
for a thorough account see Moerdijk-Reyes \cit!{3}, for a discussion 
from the point of view of functional analysis see \cit!{1}. In 
\cit!{1}, 6.6 one finds a method to recognize $C^\infty$-algebras 
among locally-m-convex algebras.

\proclaim{\nmb0{2}. Theorem}
Let $A$ be a $C^\infty$-algebra. Then $A$ is the algebra of smooth 
functions on some finite dimensional paracompact Hausdorff second 
countable manifold $M$ if and only if the following conditions are 
satisfied:
\roster
\item $A$ is \idx{\it point determined} (\cit!{3}, 4.1), so $A$ can 
       be embedded as algebra into a power $\prod_{x\in X}\Bbb R$ of 
       copies of $\Bbb R$. Equivalently the intersection of all 
       ideals of codimension 1 in $A$ is 0.
\item $A$ is \idx{\it finitely generated}, so 
       $A=C^\infty(\Bbb R^n)/I$ for some ideal 
       $I\subset C^\infty(\Bbb R^n)$.
\item For each ideal $\frak m_x$ of codimension 1 in $A$ the 
       localization $A_{\frak m_x}$ is isomorphic to the 
       $C^\infty$-algebra $C^\infty_0(\Bbb R^m)$ consisting of all 
       germs at 0 of smooth functions on $\Bbb R^m$, for some $m$.
\endroster
\endproclaim

\demo{Proof}
By condition \therosteritem2 $A$ is finitely generated, 
$A=C^\infty(\Bbb R^n)/I$; so by \cit!{3}, 4.2 
the $C^\infty$-algebra $A$ is point 
determined \therosteritem1 if and only if the ideal $I$ 
has the following property: 
$$
\text{For }f\in C^\infty(\Bbb R^n),\quad f|Z(I)=0 \text{ implies 
}f\in I,
\tag4$$ 
where $Z(I)=\bigcap \{f\i(0):f\in I\}\subset \Bbb R^n$.
Let us denote by $\{\frak m_x: x\in M\}$ the set of all ideals 
$\frak m_x$ of codimension 1 in $A$. Then $A/\frak m_x\cong \Bbb R$ 
and we write $a(x)$ for the projection of $a\in A$ in $A/\frak m_x$.
In particular we indentify the elements of $A$ with functions on $M$.
Let $c_1,\dots,c_n\in A$ by a set of generators. Then we may view 
$c=(c_1,\dots,c_n):M\to \Bbb R^n$ as a mapping such that the pullback 
$c^*(f)= f\o c = A(f)(c)$ is the quotient mapping 
$C^\infty(\Bbb R^n)\to C^\infty(\Bbb R^n)/I=A$. By condition 
\therosteritem1 $c:M\to \Bbb R^n$ is injective, and the image $c(M)$ 
equals $Z(I)=\bigcap\{f\i(0):f\in M\}$, by \thetag4. In particular, 
$c(M)$ is closed. The initial topology on $M$ with respect to all 
functions in $A$ coincides with the subspace topology induced via the 
embedding $c:M\to \Bbb R^n$, so this topology is metrizable and 
locally compact. 

Let us fix a `point' $x\in M$. The codimension 1 ideal $\frak m_x$ is a 
prime ideal, so the subset $A\setminus \frak m_x\subset A$ is closed 
under multiplication and without divisors of 0, thus
the localization $A_{\frak m_x}$ may be viewed as
the set of fractions $\frac ab$ with 
$a\in A$, $b\in A\setminus \frak m_x$; it is a local algebra with 
maximal ideal 
$\tilde{\frak m}_x = \{\frac ab:a\in\frak m_x, b\in A\setminus \frak m_x\}$.
Note that 
$\tilde\frak m_x/\tilde\frak m_x^2\cong T_0^*\Bbb R^m=\Bbb R^m$ by 
condition \thetag3.  
Now choose $a_1,\dots,a_m\in \frak m_x$ such that 
$\frac{a_1}1,\dots,\frac{a_m}1\in A_{\frak m_x}$ form a basis of 
$\tilde\frak m_x/\tilde\frak m_x^2=\Bbb R^m$, and choose 
$g_1,\dots,g_m\in C^\infty(\Bbb R^n)$ with $c^*(g_i)=a_i$. Then 
$g_i(c(x))=0$, so $g_i$ is in the codimension 1 ideal 
$\frak m_{c(x)}=\{f\in C^\infty(\Bbb R^n): f(c(x))=0\}$.
Since $c^*:C^\infty(\Bbb R^n)\to A$ induces in turn homomorphisms
$$\gather
C^\infty_{c(x)}(\Bbb R^n) = C^\infty(\Bbb R^n)_{\frak m_{c(x)}} 
\to A_{\frak m_x}\\
\Bbb R^n = T^*_{c(x)}\Bbb R^n = 
\tilde\frak m_{c(x)}/\tilde\frak m_{c(x)}^2 
\to \tilde\frak m_x/\tilde\frak m_x^2=\Bbb R^m\\
\endgather$$
and since $\frak m_{c(x)} \cong \frak m_x \oplus I$ as vector spaces, 
we may find 
functions $g_{m+1},\dots,g_n\in I$ such that the quotients 
$\frac{g_1}1,\dots,\frac{g_n}1\in C^\infty_{c(x)}(\Bbb R^n)$ map to a 
basis of $\tilde\frak m_{c(x)}/\tilde\frak m_{c(x)}^2 
=T^*_{c(x)}\Bbb R^n$. By the implicit function theorem on $\Bbb R^n$ 
the functions $g_{m+1},\dots,g_n$ are near $c(x)$ an equation of 
maximal rank for $c(M)=Z(I)$, and the functions $g_1,\dots,g_m$ 
restrict to smooth coordinates near $c(x)$ on the closed submanifold 
$c(M)=Z(I)$ of $\Bbb R^n$, and the number $m$ turns out to be a 
locally constant function on $M$. Also the functions $a_1,\dots,a_m$ 
restrict to smooth coordinates near $x$ of $M$.
\qed\enddemo

\enddocument